\documentclass[11pt,leqno,twoside]{article}

\usepackage{amsfonts,amsmath,amsthm,amssymb}
\usepackage{enumerate}
\usepackage{graphics,graphicx,subfigure}

\usepackage{color}
\usepackage{todonotes}
\usepackage{cancel}
\usepackage{url}
\usepackage{hyperref}
\usepackage{makeidx}
\usepackage{showidx}
\usepackage{multicol}        
\usepackage{xspace}
\usepackage{stmaryrd}        
\usepackage{pifont}          
\usepackage{fancybox}        
\usepackage{bm}

 \usepackage{fancyhdr}
\theoremstyle{plain}
 \theoremstyle{definition}
 \newtheorem{lem}{Lemma}
 \newtheorem{defn}[lem]{Definition}
 \newtheorem{thm}[lem]{Theorem}
 \newtheorem{prop}[lem]{Proposition}
 \newtheorem{cor}[lem]{Corollary}
 \newtheorem{notn}[lem]{Notations}
 \newtheorem{pb}[lem]{Problem}
 \newtheorem{form}[lem]{Formulae}
 
 \newtheorem*{rk}{Remark}
 \newtheorem*{com}{Comment}
 \newtheorem*{ex}{Example}
 \theoremstyle{remark}

 \newcommand{\blem}{\begin{lem}}
 \newcommand{\elem}{\end{lem}}
 \newcommand{\bdefn}{\begin{defn}}
 \newcommand{\edefn}{\end{defn}}
 \newcommand{\bthm}{\begin{thm} }
 \newcommand{\ethm}{\end{thm}}
 \newcommand{\bprop}{\begin{prop}}
 \newcommand{\eprop}{\end{prop}}
 \newcommand{\bcor}{\begin{cor}}
 \newcommand{\ecor}{\end{cor}}
 \newcommand{\bnotn}{\begin{notn}}
 \newcommand{\enotn}{\end{notn}}
 \newcommand{\bpb}{\begin{pb}}
 \newcommand{\epb}{\end{pb}}
 \newcommand{\bform}{\begin{form}}
 \newcommand{\eform}{\end{form}}
 \newcommand{\brk}{\begin{rk}}
 \newcommand{\erk}{\end{rk}}
 \newcommand{\bcom}{\begin{com}}
 \newcommand{\ecom}{\end{com}}
 \newcommand{\bex}{\begin{ex}}
 \newcommand{\eex}{\end{ex}}
 \newcommand{\bpf}{\begin{proof}}
 \newcommand{\epf}{\end{proof}}

\newcommand{\cA}{\mathcal{A}}

\newcommand{\cC}{\mathcal{C}}

\newcommand{\cE}{\mathcal{E}}

\newcommand{\cH}{\mathcal{H}}

\newcommand{\cK}{\mathcal{K}}

\newcommand{\cM}{\mathcal{M}}

\newcommand{\cP}{\mathcal{P}}

\newcommand{\cS}{\mathcal{S}}

\newcommand{\bC}{\mathbb{C}}

\newcommand{\bK}{\mathbb{K}}

\newcommand{\bN}{\mathbb{N}}

\newcommand{\bR}{\mathbb{R}}

\newcommand{\be}{\begin{equation}}
\newcommand{\ee}{\end{equation}}
\newcommand{\bal}{\begin{align}}
\newcommand{\eal}{\end{align}}
\newcommand{\ba}{\begin{align*}}
\newcommand{\ea}{\end{align*}}
\newcommand{\bmx}{\begin{matrix}}
\newcommand{\emx}{\end{matrix}}
\newcommand{\bbmx}{\begin{bmatrix}}
\newcommand{\ebmx}{\end{bmatrix}}
\newcommand{\bpmx}{\begin{pmatrix}}
\newcommand{\epmx}{\end{pmatrix}}
\newcommand{\bvmx}{\begin{vmatrix}}
\newcommand{\evmx}{\end{vmatrix}}

\newcommand{\wh}{\widehat}
\newcommand{\wt}{\widetilde}
\newcommand{\f}{\frac}

\newcommand{\inc}{\subseteq}

\newcommand{\setm}{\setminus}

\newcommand{\tr}{\mathrm{tr}}

\newcommand{\argmin}{{\rm argmin}\,}

\newcommand{\minimize}[1]{\underset{#1}{\rm minimize}\,}
\newcommand{\maximize}{{\rm maximize}\,}

\newcommand{\la}{\lambda}

\newcommand{\eps}{\varepsilon}
  
\newcommand{\TODO}[1]{{\color{red}\tiny TODO: {#1}}}     

\title{\vspace{-20mm}
Instances of Computational Optimal Recovery:\\
Dealing with Observation Errors
\medskip\hrule height 1.2pt \vspace{-6mm}}
\author{Mahmood Ettehad and Simon Foucart\footnote{S. F. is partially supported by NSF grants DMS-1622134 and DMS-1664803, and also acknowledges the NSF grant CCF-1934904.}\\ 
Texas A\&M University}
\date{\vspace{-6mm}\rule{100mm}{0.8pt}}

\newcommand\shorttitle{Instances of Computational Optimal Recovery:
Dealing with Observation Errors}
\newcommand\authors{M. Ettehad, S. Foucart}

\begin{document}
\maketitle

\vspace{-15mm}
\begin{abstract}
When attempting to recover functions from observational data,
one naturally seeks to do so in an optimal manner with respect to some modeling assumption.
With a focus put on the worst-case setting, this is the standard goal of Optimal Recovery.
The distinctive twists here are the consideration of inaccurate data through some  boundedness models and the emphasis on computational realizability.
Several scenarios are unraveled through the efficient constructions of optimal recovery maps:
local optimality under linearly or semidefinitely describable models,
global optimality for the estimation of linear functionals under approximability models,
and global near-optimality under approximability models in the space of continuous functions. 
\end{abstract}



\vspace{-5mm}
\begin{center}
\rule{100mm}{0.8pt}
\end{center}


\section{Introduction}

The investigations conducted in this article fit in the classical setting of Optimal Recovery \cite{OR}:
given observational data and {\sl a~priori} information about a function,
one attempts to approximate, in a worst-case setting, the whole function or merely a dependant quantity.
Our emphasis here is put not only on inaccurate data,
but also on computability of the approximation procedure.

Formally, functions $f$ are viewed as elements of a normed space $X$.
An educated belief about the behavior of $f$ translates into the statement that $f$ belong to a model set $\cK \inc X $ --- this is the {\sl a~priori} information.
The observational data typically take the form of evaluations of $f$ at certain points $x_1,\ldots,x_m$,
i.e., one has access to $y_1 = f(x_1),\ldots,y_m = f(x_m)$.  
In a more general and realistic situation, 
one has access to
\be
y_i = \ell_i(f) + e_i,
\qquad i \in [1:m],
\ee
where $\ell_1,\ldots,\ell_m \in X^*$ are (known) linear functionals
and the observation process is corrupted by (unknown) errors $e_1,\ldots,e_m \in \bR$.
The monograph \cite{Pla} contains substantial information about this framework under different models for the vector $e \in \bR^m$ of observation errors.  
Here, it is not viewed as random noise,
but rather assumed to belong to an uncertainty set $\cE \inc \bR^m$.
We use the terminology observation operator to denote the linear map
$L: f \in X \mapsto [\ell_1(f); \ldots; \ell_m(f)] \in \bR^m$.
This operator is considered a fixed entity,
i.e., the user does not have the leisure to select favorable observation functionals $\ell_1,\ldots,\ell_m \in X^*$.

We assume throughout the article that $Q: X \to Z$ is a linear map 
(think of $Q$ being a linear functional, in which case $Z = \bR$,
or $Q$ being the identity, in which case $Z=X$).
Our goal is to approximate the quantity of interest $Q(f)$
using only the data $y \in \bR^m$,
i.e., to produce a recovery map $R: \bR^m \to Z$ yielding a small error $\|Q(f)-R(y)\|_Z$.
With the worst-case setting in mind,
one defines two types of error, namely \vspace{-3mm}
\begin{itemize}
\item the local error of $R$ at $y$ over $\cK$ and $\cE$ is
\be
e^{\rm loc}_{\cK,\cE,Q}(L,R(y)) := \sup_{\substack{f \in \cK, \, e \in \cE\\L(f)+e=y}} \|Q(f) - R(y)\|_Z;
\ee
\item the (global) error of $R$ over $\cK$ and $\cE$ is
\be
\label{BlobErr}
e_{\cK,\cE,Q}(L,R) := \sup_{\substack{f \in \cK, \, e \in \cE}} \|Q(f) - R(L(f)+e)\|_Z.
\ee
\end{itemize}
One of the primary concerns in Optimal Recovery is to grasp how small these worst-case errors can possibly be.
This is quantified e.g. via the so-called intrinsic error 
\be
\label{IntrE}
e^*_{\cK,\cE,Q}(L) := \inf_{R: \bR^m \to Z} e_{\cK,\cE,Q}(L,R).
\ee
In this article, we put an extra emphasis on the practical computation of optimal recovery maps $R^{\rm opt}: \bR^m \to Z$,
be they \vspace{-3mm}
\begin{itemize}
\item locally optimal at $y \in L(\cK) + \cE$, in the sense that 
\be
e^{\rm loc}_{\cK,\cE,Q}(L,R^{\rm opt}(y))
= \inf_{R: \bR^m \to Z}
e^{\rm loc}_{\cK,\cE,Q}(L,R(y));
\ee
\item (globally) optimal, in the sense that
\be
e_{\cK,\cE,Q}(L,R^{\rm opt}) =
\inf_{R: \bR^m \to Z} e_{\cK,\cE,Q}(L,R).
\ee
\end{itemize}
Obviously, if a recovery map is locally optimal at any $y \in L(\cK) + \cE$,
then it is also globally optimal.
Given $y \in L(\cK) + \cE$,
we introduce the set of $f \in \cK$ consistent with the data $y$ 
and denote it as
\be
\cK_\cE(y) := \{ f \in \cK: y = L(f)+e \mbox{ for some }e \in \cE \}.
\ee
As is well-known (and not difficult to observe),
any locally optimal recovery map outputs
a so-called Chebyshev center\footnote{Note that the common definition of Chebyshev center of a set $\cS$ used here (center of the smallest ball containing~$\cS$)
differs from the one used in \cite[p.148]{BoyVan} (center of the largest ball contained in $\cS$).} of $Q(\cK_\cE(y))$, 
i.e., a center of a ball of smallest radius containing $Q(\cK_\cE(y))$.
This almost tautological observation, however,
is not enough to yield computable optimal recovery maps.
Their efficient construction constitutes one of the main points of this paper, which are listed below.\vspace{-5mm}
\begin{enumerate}
\item Examples of computable Chebyshev centers beyond the Hilbert case
(which was solved by \cite{BCDDPW} in the accurate setting under approximability models):
the results are presented in Section \ref{SecChe}
and seem to be new even in the absence of observation errors.
\item Description of globally optimal recovery maps when the quantity of interest $Q$ is a linear functional:
this extension to the inaccurate setting of a result from \cite{DFPW} appears in Section \ref{SecEstLF}
and the computational procedure proposed there is new even in the absence of observation errors when the observation functionals are not point evaluations. 
\item Construction of globally near-optimal maps for full recovery in $X = \cC(D)$:
this again extends a result from \cite{DFPW} to the inaccurate setting,
with a notable difference concerning genuinely optimal maps, see Section \ref{SecFullRec}.\vspace{-5mm}
\end{enumerate}
The remaining Section \ref{SecReminders} serves as a reminder of known facts about globally optimal recovery maps.
All our theoretical results are computationally exemplified in the reproducible {\sc matlab} file accompanying this article,
which is downloadable from the authors' webpages.

\section{Computation of Chebyshev centers}
\label{SecChe}

In this section, we uncover situations where,
using techniques from Robust Optimization~\cite{RO},
it is possible to exactly compute the Chebyshev center of the set $Q(\cK_\cE(y))$
defined by the observational data $y \in \bR^m$.
It is assumed here that
the quantity of interest $Q$ takes values in $Z = \ell_\infty^K$ 
and we write
\be
\label{FormQ}
Q: f \in X \mapsto [Q_1(f);\ldots;Q_K(f)] \in \ell_\infty^K.
\ee
The case of linear functionals is included as the special instance $K=1$.
The starting point is the observation that Chebyshev centers and radii are solutions to the 
formal optimization problem
\be
\label{HaveTo}
\minimize{z \in \bR^K, r \in \bR} \; \; r
\qquad \mbox{subject to }
\|Q(f)-z\|_\infty \le r
\quad \mbox{for all } f \in \cK_\cE(y).
\ee
We consider in this section an uncertainty set given by 
\be
\label{UncSetInf}
\cE = \cE_{\infty,\eta} 
:= \{ e \in \bR^m: \|e\|_\infty \le \eta \},
\ee
so that the condition $f \in \cK_\cE(y)$ translates into $f \in \cK$ and $y - \eta \le L(f) \le y + \eta$,
where the inequalities are understood componentwise.
As for the model set $\cK$ itself,
the method presented below essentially relies on a linear or semidefinite description for it.
This point is clarified in the following two subsections,
which are illustrated in the reproducible file.

\subsection{The model set is a polytope}

We suppose here that the set $\cK$ is a polytope in $X = \bR^n$ described in inequality form as
\be
\label{PolyMod}
\cK = \{f \in \bR^n: A f \le b \}
\quad \mbox{for some } A\in \bR^{N \times n}
\; \mbox{ and } \; b \in \bR^N.
\ee
In order to state the result of this subsection,
we define an auxiliary matrix $\wt{A} \in  \bR^{(N+2m)\times n}$ and an auxiliary vector $\wt{b} \in \bR^{N+2m}$ by
\be
\wt{A} := \bbmx A \\ \hline L \\ \hline -L \ebmx
\qquad \mbox{and} \qquad 
\wt{b} := \bbmx b \\ \hline +y + \eta \\ \hline -y + \eta \ebmx .
\ee
We also introduce the vectors $q_k \in \bR^n$, $k \in [1:K]$, that satisfy 
\be
\langle q_k,f \rangle = Q_k(f) 
\qquad \mbox{for all } f \in \bR^n.
\ee

\bthm
For the model set $\cK$ and the uncertainty $\cE$ set given in \eqref{PolyMod} and \eqref{UncSetInf},
a locally optimal recovery map over $\cK$ and $\cE$ for the quantity of interest \eqref{FormQ} outputs,
for any input $y \in L(\cK) + \cE$,
a vector $z \in \bR^K$ which is solution to the linear program
\be
\label{LocOptLP}
\minimize{\substack{z \in \bR^K, r \in \bR\\ x^{\pm,1},\ldots,x^{\pm,K} \in \bR^{N+2m}}} \; r
\qquad \mbox{subject to }
\wt{A}^\top x^{\pm,k} = \pm q_k,
\quad x^{\pm,k} \ge 0,
\quad \langle \wt{b}, x^{\pm,k} \rangle \le r \pm z_k.
\ee
\ethm

\bpf
Fixing $y \in L(\cK) + \cE$,
we first remark that the set $\cK_\cE(y)$ is a polytope described in inequality form as
\be
\cK_\cE(y) = \{ f \in \bR^n: \wt{A} f \le \wt{b} \}.
\ee
It follows that the constraint in \eqref{HaveTo} reads $\max_{f \in \bR^n} \{ \|Q(f)-z\|_\infty: \wt{A} f \le \wt{b} \} \le r$,
which in fact consists of the $2K$ linear constraints
\begin{align}
\max_{f \in \bR^n} \{ +(\langle q_k, f \rangle-z_k): \wt{A} f \le \wt{b} \} \le r,
\qquad k \in [1:K],\\
\max_{f \in \bR^n} \{ -(\langle q_k, f \rangle -z_k): \wt{A} f \le \wt{b} \} \le r,
\qquad k \in [1:K].
\end{align}
Invoking duality in linear programming
(see e.g. \cite[p.224]{BoyVan} read from the bottom up)
to transform all these $\max$-constraints
into $\min$-constraints,
the constraint in \eqref{HaveTo}
reduces to the $2K$ linear constraints
\begin{align}
\min_{x \in \bR^{N+2m}} \{ \langle \wt{b},x \rangle: \wt{A}^\top x = +q_k, \; x \ge 0 \} \le r + z_k,   
\qquad k \in [1:K],\\
\min_{x \in \bR^{N+2m}} \{ \langle \wt{b},x \rangle: \wt{A}^\top x = -q_k, \; x \ge 0 \} \le r - z_k,    
\qquad k \in [1:K].
\end{align}
Thus, the constraint in \eqref{HaveTo} is equivalent to the existence of $x^{+,1},\ldots, x^{+,K},x^{-,1},\ldots,x^{-,K} \in \bR^{N+2m}$ such that 
$\wt{A}^\top x^{\pm,k} = \pm q_k$,
$x^{\pm,k} \ge 0$, and
$\langle \wt{b}, x^{\pm,k} \rangle \le r \pm z_k$
for all $k \in [1:K]$.
Incorporating these variables into the minimization \eqref{HaveTo} leads to the announced optimization program.
We note in passing that this linear program features $K(2N+4m+1)+1$ variables,
$2Kn$ equality constraints,
and $2K(N+2m+1)$ inequality constraints.
It is therefore efficiently solvable in practice.
\epf

\subsection{The model set is the unit ball of a polynomial space}

We suppose here that the set $\cK$ is the unit ball in the space $X=\cP_n$ of algebraic polynomials of degree~$<n$ equipped with the supremum norm on $[-1,1]$. 
In other words,
\be
\label{PolynMod}
\cK = \{
f \in \cP_n: \|f\|_{\cC[-1,1]} \le 1
\}.
\ee
In order to state the result of this subsection,
we introduce the notation ${\rm Toep}(x)$ for the symmetric Toeplitz matrix built from a vector $x \in \bR^d$,
i.e.,
\be
{\rm Toep}(x) := \bbmx
x_1 & x_2 & \cdots & \cdots & x_d\\
x_2 & x_1 & x_2 & \ddots & \vdots \\
\vdots & \ddots& \ddots & \ddots & \vdots \\
\vdots & \ddots & x_2 & x_1 & x_2\\
 x_d & \cdots & \cdots & x_2 & x_1
\ebmx.
\ee
With $T_j$ denoting the $j$th degree Chebyshev polynomial of the first kind,
we now introduce the auxiliary matrices $C_1,\ldots,C_K \in \bR^{n \times n}$ and $A_1,\ldots, A_m \in \bR^{n \times n}$ defined by
\be
C_k = {\rm Toep}[Q_k(T_0); Q_k(T_1); \ldots; Q_k(T_{n-1})]
\quad \mbox{and} \quad
A_i = {\rm Toep}[\ell_i(T_0); \ell_i(T_1);\ldots; \ell_i(T_{n-1})].
\ee
The proof of the theorem below makes use of the following semidefinite duality statement,
which is a somewhat tedious application of \cite[p.452-454]{RO}.

\blem
\label{LemDual}
Given symmetric matrices $A_1,\ldots,A_m$, $B_1,\ldots,B_n$, and $C$,
the dual to the semidefinite program 
\begin{align}
\underset{M,P}{\maximize}
\; \tr[C(P-M)]
& & \mbox{subject to } &
\; \tr[B_j(P+M)] = \beta_j,
\; M,P \succeq 0,\\
\nonumber
& & \mbox{and } & 
\; \tr[A_i(P-M)] \le \alpha_i
\end{align}
is the semidefinite program
\begin{align}
\minimize{x \in \bR^n, u \in \bR^m}
\; \langle \beta, x \rangle + \langle \alpha, u \rangle
 & & \mbox{subject to } & 
\; \sum_{j=1}^n x_j B_j \succeq + \Big( C - \sum_{i=1}^m u_i A_i \Big),\\
\nonumber
& & \mbox{and } & 
\; \sum_{j=1}^n x_j B_j \succeq - \Big(C - \sum_{i=1}^m u_i A_i \Big),\\
\nonumber
& & \mbox{and } & \; u \ge 0.
\end{align}
\elem

With this lemma at hand, we can now state and prove the awaited main result of this subsection.

\bthm
For the model set $\cK$ and the uncertainty set $\cE$ given in \eqref{PolynMod} and \eqref{UncSetInf},
a locally optimal recovery map over $\cK$ and $\cE$ for the quantity of interest \eqref{FormQ} outputs,
for any input $y \in L(\cK) + \cE$,
a vector $z \in \bR^K$ which is solution to the semidefinite program
\begin{align}
\label{LocOPtSDP}
\minimize{\substack{z \in \bR^K, r \in \bR\\
x^{\pm,k} \in \bR^{n}, u^{\pm,k},v^{\pm,k} \in \bR^m}} r 
& & \mbox{s.to } & 
\; x_1^{\pm,k} + \langle u^{\pm,k}, y +\eta \rangle
-\langle v^{\pm,k}, y - \eta \rangle \le r \pm z_k,
\;  u^{\pm,k} \ge 0, \, v^{\pm,k} \ge 0,\\
\nonumber
& & \mbox{and } & 
\, {\rm Toep}(x^{\pm,k}) \succeq + \Big( \pm C_k + \sum_{i=1}^m (v^{\pm,k}_i-u^{\pm,k}_i) A_i \Big),\\
\nonumber
& & \mbox{and } & 
\, {\rm Toep}(x^{\pm,k}) \succeq - \Big( \pm C_k + \sum_{i=1}^m (v^{\pm,k}_i-u^{\pm,k}_i) A_i \Big).
\end{align}
\ethm

\bpf
It was observed in \cite[Subsection 5.3]{MinProjExt},
following ideas formulated in \cite{FouPow},
that the unit ball in $\cP_n$ admits the  semidefinite description
\begin{align}
\cK = \bigg\{
\sum_{j=0}^{n-1} \tr [D_j(P-M)] T_j
& \mbox{ for some positive semidefinite matrices }M,P \in \bR^{n \times n} \\
\nonumber
& \mbox{ that satisfy }
 \tr[D_j(P+M)]
= \delta_{0,j}
\bigg\},
\end{align}
where, for $j \in [0:n-1]$,
the symmetric matrix
\be
D_j := \bbmx
0 & \cdots & 0 & 1  & 0 & \cdots & 0\\
\vdots & 0 & \ddots & 0 & 1 & \ddots & \vdots \\
0 & \ddots & \ddots & \ddots & \ddots & \ddots & 0\\
1 & 0 & \ddots & \ddots & \ddots & 0 & 1\\
0 & 1 & \ddots & \ddots & \ddots & \ddots & 0\\
\vdots & \ddots & \ddots & 0 & \ddots & 0 & \vdots\\
0 & \cdots & 0 & 1 & 0 & \cdots & 0
\ebmx
\in \bR^{n \times n}
\ee
has $1$'s on the $j$th subdiagonal and superdiagonal and $0$'s elsewhere --- in particular $D_0$ is the $n \times n$ identity matrix.
Thus,
fixing $y \in L(\cK) + \cE$,
the set $Q(\cK_\cE(y))$ admits the semidefinite description
\begin{align}
Q(\cK_\cE(y)) = \bigg\{
\sum_{j=0}^{n-1} \tr [D_j(P-M)] Q(T_j) : & 
\; \tr [D_j(P + M)] = \delta_{0,j},
\; M,P \succeq 0,\\
\nonumber
& \; y-\eta \le \sum_{j=0}^{n-1} \tr [D_j(P-M)] L(T_j) \le y + \eta
\bigg\}.
\end{align}
Since the matrices $C_1,\ldots,C_K$ and $A_1,\ldots,A_m$ are equivalently written
as 
\be
C_k = \sum_{j=0}^{n-1} Q_k(T_j) D_j
\qquad \mbox{and} \qquad
A_i = \sum_{j=0}^{n-1} \ell_i(T_j) D_j,
\ee
we see that the constraint in \eqref{HaveTo} can be reformulated
as the $2K$ semidefinite constraints (indexed by $+$ and $-$ for each $k \in [1:K]$)
\begin{align}
\max_{M,P \in \bR^{n \times n}}
\bigg\{
\pm \tr[ C_k (P-M)] : & 
\; \tr[D_j (P + M)]
= \delta_{0,j},
\; M, P \succeq 0,\\
\nonumber
& \; y_i - \eta \le \tr [A_i (P - M)] \le y_i + \eta
\bigg\} 
\le r \pm z_k.
\end{align}
Relying on Lemma \ref{LemDual}
to transform these max-constraints into min-constraints,
the constraint in \eqref{HaveTo} reduces to the $2K$ semidefinite constraints
(indexed by $+$ and $-$ for each $k \in [1:K]$)
\begin{align}
\min_{\substack{x \in \bR^n\\ u,v \in \bR^m}}
\bigg\{ x_1 + \langle u,y+\eta \rangle - \langle v, y - \eta \rangle : 
& \; \sum_{j=0}^{n-1} x_j D_j 
\succeq + \Big( \pm C_k + \sum_{i=1}^m (v_i-u_i) A_i \Big),\\
\nonumber
& \; \sum_{j=0}^{n-1} x_j D_j 
\succeq - \Big( \pm C_k + \sum_{i=1}^m (v_i-u_i) A_i \Big),\\
\nonumber
& \; u \ge 0, v \ge 0
\bigg\} 
 \le r \pm z_k.
\end{align}
For each of these constraints,
we create extra variables $x^{\pm,k} \in \bR^{n}$, $u^{\pm,k} \in \bR^m$, and $v^{\pm,k} \in \bR^m$
to be incorporated in the optimization program \eqref{HaveTo}, which is then reformulated as
\begin{align}
\minimize{\substack{z \in \bR^K, r \in \bR\\
x^{\pm,k} \in \bR^{n}, u^{\pm,k},v^{\pm,k} \in \bR^m}} r 
& & \mbox{s.to } & 
\, x_1^{\pm,k} + \langle u^{\pm,k}, y +\eta \rangle
-\langle v^{\pm,k}, y - \eta \rangle \le r \pm z_k,
\; u^{\pm,k} \ge 0, \, v^{\pm,k} \ge 0,\\
\nonumber
& & \mbox{and } &
\, \sum_{j=0}^{n-1} x^{\pm,k}_j D_j \succeq + \Big( \pm C_k + \sum_{i=1}^m (v^{\pm,k}_i-u^{\pm,k}_i) A_i \Big),\\
\nonumber
& & \mbox{and } &
\, \sum_{j=0}^{n-1} x^{\pm,k}_j D_j \succeq - \Big( \pm C_k + \sum_{i=1}^m (v^{\pm,k}_i-u^{\pm,k}_i) A_i \Big).
\end{align}
This is indeed the announced semidefinite program,
which is solvable in practice.
\epf

\section{Global optimality over approximability models}
\label{SecReminders}

This section recollects some known ingredients that are needed later.
As such, it does not contain any new result.

\subsection{Formal reduction to the accurate setting}
\label{SecRed}

Traditional Optimal Recovery often disregards observation errors and works in the setting $e=0$.
This is because observation errors can be absorbed into the accurate setting, at least formally.
Let us recall the implicit argument (found e.g. in \cite{OR}), 
which is valid for arbitrary model and uncertainty sets $\cK$ and $\cE$.
It consists of the remark that the global error of a recovery map $R: \bR^m \to Z$ over $\cK$ and $\cE$ can be written as
\be
\sup_{\substack{f \in \cK\\e \in \cE}} \|Q(f) - R(L(f)+e)\|_Z
= \sup_{(f,e) \in \cK \times \cE} \|\wt{Q}((f,e)) - R(\wt{L}((f,e)) \|_Z,
\ee
where the quantity of interest $\wt{Q}$ 
and observation operator $\wt{L}$ are defined on the augmented space $\wt{X} = X \times \bR^m$ by
\begin{align}
\wt{Q}: & (f,e) \in \wt{X} \mapsto Q(f) \in Z,\\ 
\wt{L}: & (f,e) \in \wt{X} \mapsto L(f) + e \in \bR^m. 
\end{align}
Thus, inaccurate optimal recovery over the model and uncertainty sets $\cK$ and $\cE$
becomes optimal recovery over the model set $\wt{\cK} = \cK \times \cE \inc \wt{X}$.
This implies, for instance, that if $Q$ is a linear functional and if $\cK$ and $\cE$ are both symmetric and convex sets, then there is an optimal recovery map which is linear.

\subsection{Approximability models}

We concentrate for the rest of this article on a certain model set $\cK$ introduced in \cite{BCDDPW}.
It is given in terms of approximability by a linear space $V \inc X$ with threshold $\eps > 0$ as
\be
\label{AppModel}
\cK = \cK_{V,\eps} 
:= \{ f \in X: {\rm dist}_X(f,V) \le \eps \}.
\ee
The unit ball of $X$,
which is often considered as a model set in traditional Optimal Recovery,
corresponds to the specific choice $V =\{0\}$ and $\eps=1$.
In turn,
any symmetric convex body can be described through \eqref{AppModel} with $V =\{0\}$ and $\eps=1$,
since such a body can be viewed as the unit ball relative to some norm (namely, to its Minkowski functional).
In the case of a general space $V$, the Optimal Recovery problem under the approximability set \eqref{AppModel}
does not make sense
when its dimension exceeds the amount of data,
so one assumes that
\be
\label{Restri}
n := \dim(V) \le m.
\ee
We now recall some results valid in the absence of observation errors, see \cite[Theorems 2.1 and 3.1]{DFPW}.\vspace{-5mm}
\begin{enumerate}[(i)]
\item \label{Thm2.1}
if $Q: X \to Z$ is a linear map,
then 
the intrinsic error over the approximability set \eqref{AppModel}
satisfies
\be
\mu_{V,Q}(L) \times \eps 
\le \inf_{R: \bR^m \to Z} \sup_{f \in \cK} \|Q(f) - R(L(f))\|_Z
\le 2 \times \mu_{V,Q}(L) \times \eps,
\ee
where the indicator of compatibility between the model and the data is defined as
\be
\mu_{V,Q}(L)
:= \sup_{f \in \ker(L)} \f{\|Q(f)\|_Z}{{\rm dist}_X(f,V)};
\ee
\item 
\label{Decouple}
if $Q: X \to \bR$ is a linear functional, then the intrinsic error over the approximability set \eqref{AppModel} decouples exactly as the product of the indicator of compatibility and the approximability threshold,
i.e.,
\be
 \inf_{R: \bR^m \to Z} \sup_{f \in \cK} \|Q(f) - R(L(f))\|_Z
= \mu_{V,Q}(L) \times \eps;
\ee
\item 
\label{Thm3.1}
if $Q: X \to \bR$ is a linear functional, then
a globally optimal recovery map over the approximability set \eqref{AppModel} is provided
by the linear functional $R^{\rm opt}: y \in \bR^m \mapsto \sum_{i=1}^m a^{\rm opt}_i y_i \in \bR$, 
where the optimal weights $a^{\rm opt} \in \bR^m$ 
are precomputed (independently of $\eps>0$) as a solution to
\be 
\minimize{a \in \bR^m} 
\bigg\| Q - \sum_{i=1}^m a_i \ell_i \bigg\|_{X^*}
\qquad \mbox{subject to }
\sum_{i=1}^m a_i \ell_i(v) = Q(v)
\mbox{ for all }v \in V.
\ee
\end{enumerate}

\section{Estimation of linear functionals under approximability models}
\label{SecEstLF}

In this section,
we assume that the quantity of interest $Q$ is a linear functional. 
We place ourselves under the approximability model \eqref{AppModel}
and continue to do so throughout the rest of the article.
From now~on,
we also assume boundedness of the observation error $e \in \bR^m$, 
and hence concentrate on the uncertainty set
\be
\label{BddModel}
\cE = \cE_{p,\eta} 
:= \{ e \in \bR^m: \|e\|_p \le \eta \}
\ee
defined by an index $p \in [1,\infty]$ and parameter $\eta >0$.
It will be convenient to write $p' \in [1,p]$ for the conjugate exponent to $p$,
i.e., for $p' = p/(p-1)$ which satisfies $1/p + 1/p' = 1$.

\subsection{Description of an optimal recovery map}

The result presented in this subsection is an extension of \eqref{Thm3.1} to the inaccurate setting.
Although a dependence on $\eps > 0$ now appears,
a pleasing feature persists: the costly computation \eqref{Main} of optimal weights is performed offline once and for all.
Thus, when new data $y \in \bR^m$ comes in,
producing the associated estimate via \eqref{ProdMain} is almost immediate.
This contrasts with procedures \eqref{LocOptLP} and \eqref{LocOPtSDP},
where producing a locally optimal estimate involved a costly minimization for every new $y \in \bR^m$ coming in.

\bthm
\label{MainThm}
If $Q: X \to \bR$ is a linear functional,
then an optimal recovery map over the model set \eqref{AppModel} and the uncertainty set \eqref{BddModel} is the linear map
\be
\label{ProdMain}
R^{\rm opt}: y \in \bR^m \mapsto \sum_{i=1}^m a^{\rm opt}_i y_i \in \bR,
\ee
where the optimal weights $a^{\rm opt} \in \bR^m$ are precomputed as a solution to
\be 
\label{Main}
\minimize{a \in \bR^m} 
\bigg\| Q - \sum_{i=1}^m a_i \ell_i \bigg\|_{X^*}
+ \f{\eta}{\eps} \|a\|_{p'}
\qquad \mbox{subject to }
\sum_{i=1}^m a_i \ell_i(v) = Q(v)
\mbox{ for all }v \in V.
\ee
\ethm

\bpf
The core explanation is that,
given the approximability set \eqref{AppModel} relative to a subspace $V$ of~$X$
and the uncertainty set \eqref{BddModel} relative to an index $p \in [1,\infty]$,
the model set $\wt{\cK} = \cK \times \cE$ itself can be interpreted as an approximability set.
For this purpose, we endow the augmented space $\wt{X} = X \times \bR^m$ with the norm
\be
\| (f,e) \|_{\wt{X}}
= \max \Big\{ \|f\|_X, \f{\eps}{\eta} \|e\|_p \Big\},
\qquad (f,e) \in \wt{X}.
\ee
From there, we notice that 
\begin{align}
(f,e) \in \wt{\cK} 
& \iff \exists v \in V: \|f-v\|_X \le \eps
\mbox{ and } \|e\|_p \le \eta\\
\nonumber
& \iff \exists (v,w) \in V \times \{0\}:
\|f - v\|_X \le \eps \mbox{ and } \f{\eps}{\eta} \|e - w\|_p \le \eps\\
\nonumber
& \iff \exists (v,w) \in V \times \{0\}: \|(f,e) - (v,w)\|_{\wt{X}} \le \eps.
\end{align}
This means that $\wt{\cK}$ reduces to the approximability set 
\be
\wt{\cK}
= \big\{ (f,e) \in \wt{X}: {\rm dist}_{\wt{X}}( (f,e) , \wt{V} ) \le \eps \big\},
\qquad \quad \wt{V} = V \times \{0\}.
\ee
From the known result \eqref{Thm3.1} about the accurate setting,
we deduce that
a globally optimal recovery map is given by $R^{\rm opt}: y \in \bR^m \mapsto \sum_{i=1}^m a^{\rm opt}_i y_i \in \bR$, 
where $a^{\rm opt} \in \bR^m$ 
is a solution to
\be 
\label{OptAlgXTilde}
\minimize{a \in \bR^m} 
\bigg\| \wt{Q} - \sum_{i=1}^m a_i \wt{\ell_i} \bigg\|_{\wt{X}^*}
\qquad \mbox{subject to }
\sum_{i=1}^m a_i \wt{\ell_i}(\wt{v}) = \wt{Q}(\wt{v})
\mbox{ for all }\wt{v} \in \wt{V}.
\ee
The constraint in \eqref{OptAlgXTilde} simply reads $\sum_{i=1}^m a_i \ell_i(v) = Q(v)$
for all $v \in V$
 because any $\wt{v} \in \wt{V}$ takes the form $\wt{v} = (v,0)$ for some $v \in V$.
As for the objective function,
it transforms into
\begin{align}
\label{ObjFunTrans}
\sup_{\|(f,e)\|_{\wt{X}} \le 1}  \bigg| \wt{Q}((f,e)) - \sum_{i=1}^m a_i \wt{\ell_i}((f,e)) \bigg|
& = \sup_{\substack{\|f\|_X \le 1\\ \|e\|_p \le \eta/\eps }}
\bigg| Q(f) - \sum_{i=1}^m a_i (\ell_i(f)+e_i) \bigg|\\
\nonumber
& = \sup_{\|f\|_X \le 1}\bigg| Q(f) - \sum_{i=1}^m a_i \ell_i(f) \bigg|
+ \sup_{\|e\|_p \le \eta / \eps} \bigg| \sum_{i=1}^m a_i e_i \bigg|\\
\nonumber
& = \bigg\| Q - \sum_{i=1}^m a_i \ell_i \bigg\|_{X^*}
+ \f{\eta}{\eps} \|a\|_{p'}.
\end{align}
The result is now fully justified by 
substituting \eqref{ObjFunTrans} as the objective function in \eqref{OptAlgXTilde}
while taking the simplified form of the constraint into account.
\epf

\brk
When $Q: X \to Z$ is an arbitrary linear map,
the interpretation of $\wt{\cK} = \cK \times \cE$ as an approximability set also
implies, by \eqref{Thm2.1}, that the intrinsic error over $\cK$ and $\cE$ satisfies 
\be
\mu_{\wt{V},\wt{Q}}(\wt{L}) \times \eps 
\le \inf_{R: \bR^m \to Z} \sup_{\substack{f \in \cK \\ e \in \cE}} \|Q(f) - R(L(f)+e)\|_Z
\le 2 \times \mu_{\wt{V},\wt{Q}}(\wt{L}) \times \eps,
\ee
where the indicator of compatibility now depends on $\eps>0$ (unless $\eta$ is proportional to $\eps$) via
\begin{align}
\label{MuDependsEps}
\mu_{\wt{V},\wt{Q}}(\wt{L}) & = \sup_{(f,e) \in \ker(\wt{L})}
\f{\|\wt{Q}((f,e))\|_Z}{{\rm dist}_{\wt{X}}((f,e),\wt{V})}
= \sup_{Lf+e = 0} \f{\|Q(f)\|_Z}{\max\{ {\rm dist}_X(f,V), \f{\eps}{\eta} \|e\|_p \} }\\
\nonumber
& = \sup_{f \in X} \f{\|Q(f)\|_Z}{\max\{ {\rm dist}_X(f,V), \f{\eps}{\eta} \|Lf\|_p \} }. 
\end{align}
This supremum over $f \in X$ is larger than or equal to the supremum over $f \in \ker(L)$,
which leads to the intuitive fact that the `inaccurate' indicator $\mu_{\wt{V},\wt{Q}}(\wt{L})$ is larger than or equal to the `accurate' indicator $\mu_{V,Q}(L)$.
It is also worth pointing out the fact that 
\be
\label{MuQMuI}
\mu_{\wt{V},\wt{Q}}(\wt{L})
\le \|Q\|_{X \to Z} \times \mu_{\wt{V},\wt{I}}(\wt{L}).
\ee
\erk

\brk
By suitably modifying the approximability set,
the result of Theorem \ref{MainThm}
can be pushed beyond the 
restriction \eqref{Restri} 
imposing some underparametrization.
For details, we refer to the companion article \cite{Fou},
which introduces and analyzes the model sets
\begin{align}
\cK & = \big\{ f \in X: \; {\rm dist}_X(f,V) \le \eps \mbox{ and } \|f\|_X \le \kappa \big\}, \\
\cK & = \big\{ f \in X: \; \exists v \in V \mbox{ with } \|f-v\|_X \le \eps \mbox{ and } \|v\|_X \le \kappa \big\}. 
\end{align}
\erk

\subsection{Computational realization for $X=\cC[-1,1]$}

Unless the dual norm of $X$ can be practically handled,
the value of Theorem \ref{MainThm} would remain at the theoretical level only.
The task of solving the optimization program \eqref{Main} is probably easiest when $X$ is a reproducing kernel Hilbert space.
We do not pursue this direction, which is really close to \cite[Subsection 5.2]{DFPW}.
Instead, we consider in this subsection
the important situation $X=\cC[-1,1]$.
We shall reveal that solving \eqref{Main} is computationally feasible in this situation, too.
Notice first that,
in the typical case emphasized in \cite{DFPW}
where the observation functionals $\ell_1,\ldots,\ell_m$ are point evaluations at distinct $x_1,\ldots,x_m \in [-1,1]$, 
the task at hand is relatively easy,
since the objective function of~\eqref{Main} reduces for a generic $Q$ to  $\|a\|_1 + (\eta/\eps) \|a\|_{p'}$,
up to the additive constant $\|Q\|_{\cC[-1,1]^*}$.
Our focus here is on observation functionals that take the general form
\be
\ell_i(f) = \int_{-1}^1 f(x) d\la_i(x),
\qquad f \in \cC[-1,1],
\ee
for some signed Borel measures $\la_1,\ldots,\la_m$ on $[-1,1]$.
As a guiding example developed in our {\sc matlab} reproducible,
and similarly to a scenario considered in \cite{AHS},
we can think of $V$ as the space $\cP^{\rm odd}_{2n}$ of odd algebraic polynomials of degree $<2n$
and of the observation functionals $\ell_1,\ldots,\ell_m$ as Fourier measurements with, say, 
$d\la_i(x) =  \sin(i \pi x) dx$.
Let us also write the linear functional $Q$ as
\be
Q(f) = \int_{-1}^1 f(x) d\rho(x),
\qquad f \in \cC[-1,1],
\ee
for some signed Borel measure $\rho$ on $[-1,1]$.
The main optimization problem \eqref{Main} then turns into
\be
\label{Main1}
\minimize{a \in \bR^m} 
\int_{-1}^1 d\bigg| \rho - \sum_{i=1}^m a_i \la_i \bigg|
+ \f{\eta}{\eps} \|a\|_{p'}
\qquad \mbox{subject to }
\sum_{i=1}^m a_i \ell_i(v_j) = Q(v_j),
\quad j \in [1:n],
\ee
where $(v_1,\ldots,v_n)$ denotes a basis for $V$.
The latter constraint reads $Ma=b$,
where the matrix $M \in \bR^{n \times m}$ and the vector $b \in \bR^n$ have entries
\be
\label{M&b}
M_{j,i} = \ell_i(v_j)
\qquad \mbox{and} \qquad
b_j = Q(v_j),
\qquad i \in  [1:m],
\quad j \in [1:n].
\ee
Let us introduce as slack variables the nonnegative Borel measures $\nu^+$ and $\nu^-$ involved in the Jordan decomposition $\nu = \nu^+ - \nu^-$ of $\nu := \rho - \sum_{i=1}^m a_i \la_i$,
so that the problem \eqref{Main1} is equivalent to
\be
\label{Main2}
\minimize{\substack{a \in \bR^m \\ \nu^+,\nu^-}} \int_{-1}^1 d(\nu^+ + \nu^-) + \f{\eta}{\eps} \|a\|_{p'}
\qquad \mbox{subject to } Ma = b,
\; \nu^+ - \nu^- = \rho - \sum_{i=1}^m a_i \la_i.
\ee
Next, replacing the measures $\nu^+$ and $\nu^-$
by their infinite sequences $z^+ = \cM_\infty(\nu^+) \in \bR^\bN$ and $z^- = \cM_\infty(\nu^-) \in \bR^\bN$ of  moments defined by
\be
z^\pm_k = 
\int_{-1}^{1} T_{k-1}(x) d\nu^{\pm}(x),
\qquad k \ge 1,
\ee
the problem \eqref{Main2} becomes equivalent\footnote{the equivalence is based on the discrete trigonometric moment problem, see \cite{FouLas} for details.} to the infinite semidefinite program
\begin{align}
\label{TrueSDP}
\minimize{\substack{a \in \bR^m \\ z^+,z^- \in \bR^\bN}} \; z^+_1 + z^-_1 + \f{\eta}{\eps} \|a\|_{p'}
& & \mbox{subject to } & \; Ma = b, \; z^+-z^- = \cM_\infty \Big(\rho-\sum_{i=1}^m a_i \la_i \Big),\\
\nonumber
& & \mbox{and } & \; {\rm Toep}_\infty(z^\pm) \succeq 0. 
\end{align}
Instead of solving this infinite optimization program,
we truncate it to a level $N$ and solve instead the resulting finite semidefinite program
\begin{align}
\label{TrunSDP}
\minimize{\substack{a \in \bR^m \\ z^+,z^- \in \bR^N}} \; z^+_1 + z^-_1 + \f{\eta}{\eps} \|a\|_{p'}
& & \mbox{subject to } & \; Ma = b, \; z^+-z^- = \cM_N \Big(\rho-\sum_{i=1}^m a_i \la_i \Big),\\
\nonumber
& & \mbox{and } & \; {\rm Toep}_N(z^\pm) \succeq 0.
\end{align}
The rest of this subsection is devoted to justifying in a quantitative way that the minimal value of this truncated problem converges to the minimal value of the original problem.
We also justify, although not quantitatively,
that the vectors $a^{(N)} \in \bR^m$ obtained by solving \eqref{TrunSDP} converge as $N \to \infty$ to a vector $a^{\rm opt} \in \bR^m$ minimizing \eqref{Main}.

\bthm
\label{ThmConv}
If the optimization program \eqref{Main} has a unique minimizer,
then the latter is the limit of any sequence $(a^{(N)})_{N \ge 1}$ 
obtained by solving \eqref{TrunSDP} for each $N \ge 1$.
Without uniqueness, it holds that
any subsequence of $(a^{(N)})_{N \ge 1}$ admits a subsequence converging to one of the minimizers of~\eqref{Main}.
\ethm

\bpf
The first part of the theorem follows from the second part:
it is indeed well-known that the convergence of a sequence to a given point is guaranteed as soon as any of its subsequences admits a subsequence converging to that point. 

To establish the second part,
let $\alpha^{(N)} \in \bR$ and  $(a^{(N)},z^{+,(N)},z^{-,(N)}) \in \bR^m \times \bR^{N} \times \bR^{N}$ denote, for each $N \ge 1$, the minimum value and some minimizer of \eqref{TrunSDP}, respectively.
We write $z^{\pm,((N))} \in \bR^\bN$ for the infinite vectors obtained by padding the finite vectors $z^{\pm,(N)} \in \bR^N$ with zeros.
Let us now consider a subsequence $\big( (a^{(N_k)},z^{+,((N_k))},z^{-,((N_k))}) \big)_{k \ge 1}$ of the whole $(\ell_\infty^m \times \ell_\infty^\bN \times \ell_\infty^\bN)$-valued sequence $\big( (a^{(N)},z^{+,((N))},z^{-,((N))})  \big)_{N \ge 1}$.
Our objective is  to show that there exist a subsequence $\big( (a^{(N_{k_\ell})},z^{+,((N_{k_\ell}))},z^{-,((N_{k_\ell}))}) \big)_{\ell \ge 1}$
and a minimizer $(\wt{a},\wt{z}^{+},\wt{z}^{-})$ of \eqref{TrueSDP} such that $a^{(N_{k_\ell})}$ converges to~$\wt{a}$ as $\ell \to \infty$.
To this end, we start by observing that the sequence $(\alpha^{(N_k)})_{k \ge 1}$ is nondecreasing and bounded  by the minimal value $\alpha^{\rm opt}$ of \eqref{TrueSDP}:
firstly, 
the inequality $\alpha^{(N_k)} \le \alpha^{(N_{k+1})}$ 
follows from the feasibility of $(a^{(N_{k+1})},z^{+,(N_{k+1})}_{[1:N_k]},z^{-,(N_{k+1})}_{[1:N_k]})$
for \eqref{TrunSDP} specified to $N=N_k$, so that 
\be
\alpha^{(N_k)}
\le z^{+,(N_{k+1})}_1 + z^{-,(N_{k+1})}_1 + \f{\eta}{\eps} \|a^{(N_{k+1})}\|_{p'} = \alpha^{(N_{k+1})};
\ee
secondly, the inequality $\alpha^{(N_k)} \le \alpha^{\rm opt}$ similarly follows from the feasibility of $(a^{\rm opt},z^{+,\rm opt}_{[1:N_k]},z^{-,\rm opt}_{[1:N_k]})$ for~\eqref{TrunSDP} specified to $N=N_k$,
where evidently $(a^{\rm opt},z^{+,\rm opt},z^{-,\rm opt})$ represents some minimizer of~\eqref{TrueSDP}.
We continue by remarking that the $\ell_\infty^m$-valued sequence $(a^{(N_k)})_{k \ge 1}$ is bounded: this is a consequence of
\be
\|a^{(N_k)}\|_\infty \le \|a^{(N_k)}\|_{p'}
\le \f{\eps}{\eta} \Big( z^{+,(N_k)}_1 + z^{-,(N_k)}_1 + \f{\eta}{\eps} \|a^{(N_k)}\|_{p'} \Big)
 = \f{\eps}{\eta} \alpha^{(N_k)} 
 \le \f{\eps}{\eta} \alpha^{\rm opt}.
\ee
We also note that the positive semidefiniteness of ${\rm Toep}_{N_k}(z^{\pm,(N_k)})$ implies that,
for any $j \in [1:N_k]$,
\be
\big| z^{\pm,(N_k)}_j \big| \le z^{\pm,(N_k)}_1 \le z^{+,(N_k)}_1 + z^{-,(N_k)}_1 + \f{\eta}{\eps} \|a^{(N_k)}\|_{p'}
= \alpha^{(N_k)} \le \alpha^{\rm opt}.
\ee
Thus,
the $\ell_\infty^\bN$-valued sequences $(z^{\pm,((N_k))})_{k \ge 1}$
are also bounded.
These last two facts guarantee 
(in particular by the sequential compactness Banach--Alaoglu theorem)
that $(a^{(N_k)})_{k \ge 1}$ admits a convergent subsequence in the standard topology of $\ell_\infty^m$
and that $(z^{\pm,((N_k))})_{k \ge 1}$ admit convergent subsequences in the weak-star topology of $\ell_\infty^\bN$.
We denote the resulting convergent subsequence and its limit by 
$\big( (a^{(N_{k_\ell})},z^{+,((N_{k_\ell}))},z^{-,((N_{k_\ell}))}) \big)_{\ell \ge 1}$ 
and $(\wt{a},\wt{z}^+,\wt{z}^-)$, respectively.
It remains to prove that the triple $(\wt{a},\wt{z}^+,\wt{z}^-)$ is a minimizer of \eqref{TrueSDP}.
Since the weak-star convergence implies that
$z^{\pm,(N_{k_\ell})}_j \to \wt{z}^{\pm}_j$ for all $j \ge 1$,
writing the constraints of \eqref{TrunSDP} specified to $N=N_{k_\ell}$ and passing to the limit as $\ell \to \infty$
shows that the triple 
is feasible for \eqref{TrueSDP}.
It is also a minimizer for this program,
by virtue of
\be
\label{IsMin}
\wt{z}^+_1 + \wt{z}^-_1 + \f{\eta}{\eps} \|\wt{a}\|_{p'} 
= \lim_{\ell \to \infty} \Big( z^{+,(N_{k_\ell})}_1 + z^{-,(N_{k_\ell})}_1 + \f{\eta}{\eps} \|a^{(N_{k_\ell})}\|_{p'} \Big)
= \lim_{\ell \to \infty} \alpha^{(N_{k_\ell})}
\le \alpha^{\rm opt}.
\ee
Our objective is now established, so
the second part of theorem is proved.
\epf

Theorem \ref{ThmConv} does not tell us how to choose $N$ in order to reach a prescribed accuracy on $\|a^{\rm opt} - a^{(N)}\|$, not even on $\alpha^{\rm opt} - \alpha^{(N)}$.
The observation below provides such a quantitative estimate,
although it is an {\sl a~prosteriori} estimate,
in the sense that a bound on $\alpha^{\rm opt} - \alpha^{(N)}$ can be evaluated only after solving \eqref{TrunSDP} for a particular $N$
--- if the accuracy is not satisfactory,
one would solve \eqref{TrunSDP} again for a higher~$N$.

\bthm
\label{ThmQuant}
For any $N \ge 1$, one has
\be
\alpha^{(N)} \le \alpha^{\rm opt}
\le \alpha^{(N)} + \delta^{(N)},
\ee
where $\delta^{(N)} \ge 0$ is a computable quantity clustering to zero defined by
\be
\delta^{(N)} := \bigg[ \Big\| Q - \sum_{i=1}^m a_i^{(N)} \ell_i \Big\|_{\cC[-1,1]^*}
+ \f{\eta}{\eps} \|a^{(N)}\|_{p'}
\bigg] - \alpha^{(N)}.
\ee
\ethm

\bpf
The leftmost inequality was already justified implicitly in the proof of Theorem \ref{ThmConv}.
For the rightmost inequality,
we simply notice that $a^{(N)}$ is feasible for \eqref{Main},
so that
\be
\alpha^{\rm opt} \le \Big\| Q - \sum_{i=1}^m a_i^{(N)} \ell_i \Big\|_{\cC[-1,1]^*}
+ \f{\eta}{\eps} \|a^{(N)}\|_{p'} 
= \alpha^{(N)} + \delta^{(N)}.
\ee
The fact that $\delta^{(N)}$ clusters to zero
follows from Theorem \ref{ThmConv} and its proof:
the term in square brackets clusters to $\alpha^{\rm opt}$ (because $a^{(N)}$ clusters to a minimizer $a^{\rm opt}$ of \eqref{Main})
and the term $\alpha^{(N)}$ also converges to $\alpha^{\rm opt}$
(because the sequence $(\alpha^{(N)})_{N \ge 1}$ is nondecreasing and bounded above by $\alpha^{\rm opt}$,
hence convergent, and its limit cannot be smaller than $\alpha^{\rm opt}$, as a consequence of \eqref{IsMin}).
\epf

\brk
Even without an estimate of $\|a^{\rm opt} - a^{(N)}\|$,
using a solution $a^{(N)} \in \bR^m$ of \eqref{TrunSDP} instead of a solution $a^{\rm opt} \in \bR^m$ of \eqref{TrueSDP}
yields a recovery map $R^{(N)}: y \in \bR^m \mapsto \sum_{i=1}^m a^{(N)}_i y_i \in \bR$ which is almost optimal.
Indeed, the worst-case error of $R^{(N)}$ over $\cK$ and $\cE$ satisfies
\begin{align}
\sup_{\substack{f \in \cK\\ e \in \cE}}
\big| Q(f) - R^{(N)}(L(f)+e) \big|
& \le \Big( \eps \, \Big\|  Q -  \sum_{i=1}^m a^{(N)}_i \ell_i \Big\|_{\cC[-1,1]^*}
+ \eta \, \|a^{(N)} \|_{p'} \Big) 
= (\alpha^{(N)} + \delta^{(N)}) \times \epsilon \\
\nonumber
& \le (\alpha^{\rm opt} + \delta^{(N)}) \times \epsilon, 
\end{align}
which is only $\delta^{(N)} \times \eps$ above the minimal worst-case error.
\erk

\section{Recovery of continuous functions under approximability models}
\label{SecFullRec}

In this section, we fix $X=\cC(D)$ for some compact domain $D$ and we consider the quantity of interest $Q=I_{\cC(D)}$,
i.e., we target the full recovery of functions $f \in \cC(D)$.
We will uncover a practical construction of linear recovery maps that are near-optimal rather than genuinely optimal.
The construction will follow closely an idea from \cite[Subsection 4.3]{DFPW}.
However, we begin by highlighting that the (unpractical) construction of  a linear genuinely optimal recovery map which was presented in \cite[Subsection~4.2]{DFPW}
does not carry over from the accurate setting to the inaccurate setting.

\subsection{Discontinuity of optimal weights}

An optimal recovery map $R^{\rm opt}: \bR \to \cC(D)$ was constructed in \cite{DFPW} as follows:
for each $x \in D$, solve the minimization problem \eqref{Main} for 
the quantity of interest $Q_x$ defined by $Q_x(f)=f(x)$,
thus producing a (carefully selected) minimizer $a^{\rm opt}(x)$;
then, with $a^{\rm opt}$ denoting the function $x \in D \mapsto a^{\rm opt}(x) \in \bR^m$,
consider the map $R^{\rm opt}$ defined for $y \in \bR^m$ by $R^{\rm opt}(y) = \sum_{i=1}^m y_i a_i^{\rm opt}$;
finally, establish the optimality of $R^{\rm opt}$ by relying on the critical fact that it takes values into $\cC(D)$. 
It is this fact that does not carry over to the inaccurate setting.
Precisely,
the function $a^{\rm opt}$ is not continuous in general,
as formalized below.

\bprop
Let $p \in (1,\infty)$, let $\eta \le \eps$, and let $\ell_1,\ldots,\ell_m$ be observation functionals that are point evaluations at distinct points $x_1,\ldots,x_m \in D$.
For $k \in [1:m]$,
as $x \in D \setm \{x_k\}$ converges to the evaluation point $x_k$,
it is not guaranteed that $a^{\rm opt}(x)$ converges to $a^{\rm opt}(x_k)$.
\eprop

\bpf
Firstly, 
when $p \in (1,\infty)$,
we note that $a^{\rm opt}(x)$ is uniquely defined for any $x \in D$ due to the strict convexity of the objective function in \eqref{Main}.
Secondly, we point out that $a^{\rm opt}(x_k)$ coincides with $e_k = [0;\ldots;0;1;0,\ldots;0]$,
i.e., that $e_k$ satisfies the appropriate constraint and minimizes the quantity
\be
\label{xxk}
\Big\| Q_{x_k} - \sum_{i=1}^m a_i \ell_i \Big\|_{\cC(D)^*}
+ \f{\eta}{\eps} \|a\|_{p'} = |1-a_k| + \sum_{i=1,i\not= k}^m |a_i| + \f{\eta}{\eps} \|a\|_{p'},
\ee
the latter being true because,
when $\eta \le \eps$,
equality occurs for $a = e_k$ in
\be
|1-a_k| + \sum_{i=1,i\not= k}^m |a_i| + \f{\eta}{\eps} \|a\|_{p'}
\ge \f{\eta}{\eps} |1-a_k| + 0+ \f{\eta}{\eps} |a_k| \ge \f{\eta}{\eps}.
\ee
Thirdly, 
we are going to prove by contradiction that 
in general $a^{\rm opt}(x) \not\to e_k$ as $x \to x_k$ with $x \not= x_k$,
keeping in mind that
$a^{\rm opt}(x)$ minimizes the quantity
\be
\label{xnotxk}
\Big\| Q_{x} - \sum_{i=1}^m a_i \ell_i \Big\|_{\cC(D)^*}
+ \f{\eta}{\eps} \|a\|_{p'} = 
1 + \|a\|_1 + \f{\eta}{\eps} \|a\|_{p'}
\ee
among all $a \in \bR^m$ satisfying $Ma = b(x)$ --- here, $b(x) \in \bR^n$ denotes the vector defined in \eqref{M&b} for the quantity of interest $Q = Q_x$.
Now let $a \in \bR^m$ satisfying $Ma = b(x_k)$ and let us consider $\wh{a} := a + M^\dagger(b(x)-b(x_k))$,
where $M^\dagger \in \bR^{m \times n}$ is a pseudoinverse of $M \in \bR^{n \times m}$.
In view of $M \wh{a} = b(x)$, we derive that
\be
\|\wh{a}\|_1 + \f{\eta}{\eps} \|\wh{a}\|_{p'}
\ge \|a^{\rm opt}(x)\|_1 + \f{\eta}{\eps} \|a^{\rm opt}(x)\|_{p'}.
\ee
Since $b(x) \to b(x_k)$ as $x \to x_k$,
if we had $a^{\rm opt}(x) \to e_k$,
then passing to the limit would give
\be
\|a\|_1 + \f{\eta}{\eps} \|a\|_{p'}
\ge 1 + \f{\eta}{\eps}.
\ee
Thus, it would hold that the minimum of $\|a\|_1 + (\eta/\eps) \|a\|_{p'}$ subject to $Ma = b(x_k)$
is always $1+(\eta/\eps)$.
But this fact is easily invalidated numerically,
see the reproducible file for the case $V = \cP_n$.
\epf

\brk
For $p \in \{1,\infty\}$,
the optimal weights may not be uniquely defined.
Consequently,
a relevant question pertains to the possibility of selecting a minimizer $a^{\rm opt}(x)$ of \eqref{Main} for $Q = Q_x$ in such a way that the resulting function $a^{\rm opt}$ is continuous.
If we insist on the intuitive selection $a^{\rm opt}(x_k) = e_k$, 
then the existence of a continuous selection implies,
as in the previous argument, the fact that $1+(\eta/\eps)$ is the minimum of $\|a\|_1 + (\eta/\eps) \|a\|_{p'}$ subject to $Ma = b(x_k)$.
This fact can be invalidated numerically for $p=1$, i.e., for $p' = \infty$.
However, for $p = \infty$, i.e., for $p' = 1$, a continuous selection does exist provided the space $V$ contains the constant functions.
This was proved in \cite[Theorem 4.2]{DFPW} in the case $\eta = 0$.
Denoting by $a^{\rm opt}$ this continuous selection,
we claim that it is also a continuous selection of minimizers of \eqref{Main} in the case $\eta > 0$ and $p'=1$.
Notice indeed that
$a^{\rm opt}(x)$ minimizes $\|a\|_1$ subject to $Ma = b(x)$ for $x \not\in \{ x_1,\ldots,x_m \}$, see \eqref{xnotxk} with $\eta = 0$,
and that $a^{\rm opt}(x_k)$ minimizes $|1-a_k| + \sum_{i\not= k}|a_i|$ subject to $Ma = b(x_k)$, see \eqref{xxk} with $\eta =0$, 
as well as $\|a\|_1$ subject to $Ma = b(x_k)$ by continuity.
It then follows that
$a^{\rm opt}(x)$ minimizes $\big\| Q_{x} - \sum_{i} a_i \ell_i \big\|_{\cC(D)^*}
+ (\eta/\eps) \|a\|_{p'}$ subject to $Ma = b(x)$ for $x \not\in \{ x_1,\ldots,x_m \}$, by virtue of \eqref{xnotxk} with $p'=1$,
and that $a^{\rm opt}(x_k)$ minimizes 
$\big\| Q_{x_k} - \sum_{i} a_i \ell_i \big\|_{\cC(D)^*}
+ (\eta/\eps) \|a\|_{p'}$ subject to $Ma = b(x_k)$, by virtue of \eqref{xxk} with $p'=1$.
In summary, the vector $a^{\rm opt}(x)$ is a minimizer of \eqref{Main} for any $x \in D$, as claimed.
\erk

\subsection{Practical construction of linear near-optimal maps}

Even though the straightforward construction of a genuinely optimal recovery map cannot be reproduced
in the inaccurate setting,
we reveal in this subsection that,
if one settles for near-optimal recovery maps,
then efficient constructions are available.
All is needed 
are linear functionals $Q_1,\ldots,Q_{\bar{n}}$ with $\|Q_j\|_{\cC(D)^*} \le 1$ and
functions $u_1,\ldots,u_{\bar{n}} \in \cC(D)$
such that the linear operator $P: \cC(D) \to \cC(D)$ defined by
\be
\label{QP}
P(f) = \sum_{j=1}^{\bar{n}} Q_j(f) u_j,
\qquad f \in \cC(D),
\ee
obeys the reproducing condition
\be
P(v) = v \qquad \mbox{ for all }v \in V,
\ee
as well as, for some $\gamma \ge 1$, the boundedness condition
\be
\|P\|_{\cC(D) \to \cC(D)}
\le 
\bigg\| \sum_{j=1}^{\bar{n}} |u_j| \bigg\|_{\cC(D)} \le \gamma.
\ee
For $D = [-1,1]$ and $V=\cP_n$,
such quasi-interpolant operators $P$ exist with ${\bar{n}} = C_\gamma n$,
with $Q_1,\ldots,Q_{\bar{n}}$ being point evaluations,
and with $u_1,\ldots,u_{\bar{n}}$ being polynomials, 
see \cite[Subsection~4.3.1]{DFPW}. 

\bthm
For the model set $\cK$ and the uncertainty set $\cE$ given in \eqref{AppModel} and \eqref{BddModel},
a near-optimal recovery map for the full approximation problem over $\cK$ and $\cE$ is provided by the linear map
\be
\label{RNear}
R^{\rm near}: y \in \bR^m \mapsto
\sum_{i=1}^m y_i a^{\rm near}_i \in \cC(D),
\qquad  \quad
a^{\rm near}_i := \sum_{j=1}^{\bar{n}} a_i^{(j)} u_j \in \cC(D),
\ee
where the vectors $a^{(j)} \in \bR^m$ are solutions to \eqref{Main} with $Q = Q_j$.
\ethm

\bpf
In view of \eqref{Thm2.1},
we aim at proving that there is a constant $C_\gamma \ge 1$ such that
\be
\label{ObjNear}
\sup_{\substack{f \in \cK\\ e \in \cE}}
\|f -R^{\rm near}(L(f)+e)\|_{\cC(D)}
\le C_\gamma \times \mu_{\wt{V},\wt{I}}(\wt{L}) \times \eps .
\ee
Let us first remark that,
for any $j \in [1:{\bar{n}}]$,
the defining property of the $a^{(j)} \in \bR^m$
yields 
\be
\label{EqMin}
\bigg\| Q_j - \sum_{i=1}^m a_i^{(j)} \ell_i \bigg\|_{\cC(D)^*} + \f{\eta}{\eps} \|a^{(j)}\|_{p'}
= \mu_{\wt{V},\wt{Q_j}}(\wt{L}),
\ee
as well as the identity
\be
\label{Reprod}
\sum_{i=1}^m a_i^{(j)} \ell_i(v) = Q_j(v)
\qquad \mbox{for all } v \in V.
\ee
The latter implies that,
for any $v \in V$,
\be
R^{\rm near}(L(v))  = \sum_{i=1}^m \ell_i(v)
\Big( \sum_{j=1}^{\bar{n}} a_i^{(j)} u_j \Big)
= 
\sum_{j=1}^{\bar{n}} \Big( \sum_{i=1}^m a_i^{(j)} \ell_i(v) \Big) u_j
= \sum_{j=1}^{\bar{n}} Q_j(v) u_j = P(v)  = v.
\ee
Let now $f \in \cK$ and $e \in \cE$ be fixed.
Given $v \in V$ such that
$h:= f-v$ satisfies
$\|h\|_{\cC(D)} \le \eps$,
we have
\be
\label{NeOp0}
f - R^{\rm near}(L(f)+e) = f-v - R^{\rm near}(L(f-v) +e)
= [h - R^{\rm near}(L(h))] - [R^{\rm near}(e)].
\ee
The second term in square brackets applied to $x \in D$ is bounded as
\begin{align}
\label{NeOp1}
\big |R^{\rm near}(e)(x) \big| 
& = \Big| \sum_{i=1}^m e_i a_i^{\rm near}(x) \Big|
\le \sum_{i=1}^m |e_i| \sum_{j=1}^{\bar{n}} |a_i^{(j)}| |u_j(x)|
= \sum_{j=1}^{\bar{n}} \Big( \sum_{i=1}^m |e_i| |a_i^{(j)}| \Big) |u_j(x)|\\
\nonumber
& \le \sum_{j=1}^{\bar{n}} \|e\|_p \|a^{(j)}\|_{p'} |u_j(x)|
\le \eta \sum_{j=1}^{\bar{n}} \|a^{(j)}\|_{p'} |u_j(x)|.
\end{align}
As for the first term in square brackets applied to $x \in D$, it is bounded as
\begin{align}
\label{NeOp2}
\big| [h  - R^{\rm near}(L(h))]& (x) \big| 
= 
 \Big| h(x) - P(h)(x) + P(h)(x) - \sum_{i=1}^m \ell_i(h) \sum_{j=1}^{\bar{n}} a_i^{(j)} u_j(x) \Big|\\
 \nonumber
& = \Big| (I-P)(h)(x) + \sum_{j=1}^{\bar{n}}
\Big( Q_j(h) - \sum_{i=1}^m a_i^{(j)} \ell_i(h) \Big) u_j(x) \Big|\\
 \nonumber
 & \le \Big| (I-P)(h)(x) \Big|
 + \sum_{j=1}^{\bar{n}} \Big| \Big( Q_j - \sum_{i=1}^m a_i^{(j)} \ell_i \Big)(h) \Big| |u_j(x)|\\
 \nonumber
 & \le
 \|I-P\|_{\cC(D) \to \cC(D)} \|h\|_{\cC(D)} + \sum_{j=1}^{\bar{n}} \Big\| Q_j - \sum_{i=1}^m a_i^{(j)} \ell_i \Big\|_{\cC(D)^*} \|h\|_{\cC(D)} |u_j(x)|\\
 \nonumber
& \le
\bigg(
1 + \gamma
+ \sum_{j=1}^{\bar{n}} \Big\| Q_j - \sum_{i=1}^m a_i^{(j)} \ell_i \Big\|_{\cC(D)^*} |u_j(x)|
\bigg) \times \eps.
\end{align}
Substituting \eqref{NeOp2} and \eqref{NeOp1} into \eqref{NeOp0} leads,
for any $x \in D$, to
\begin{align}
\big| [f - R^{\rm near} (L(f)+e)](x) \big| &  \le 
\bigg(
1+\gamma
+ \sum_{j=1}^{\bar{n}} \Big( \Big\| Q_j - \sum_{i=1}^m a_i^{(j)} \ell_i \Big\|_{\cC(D)^*} + \f{\eta}{\eps} \|a^{(j)}\|_{p'} \Big) |u_j(x)|
\bigg) \times \eps\\
\nonumber
& =
\bigg(
1+\gamma
+ \sum_{j=1}^{\bar{n}} \mu_{\wt{V},\wt{Q_j}}(\wt{L}) |u_j(x)|
\bigg) \times \eps \le \bigg( 1+\gamma + \gamma \mu_{\wt{V},\wt{I}}(\wt{L}) \bigg) \times \eps,
\end{align}
where we have used \eqref{MuQMuI}
for $Q = Q_j$ in the last step.
In view of $\mu_{\wt{V},\wt{I}}(\wt{L}) \ge 1$,
taking the supremum over $x \in D$
and then over $f \in \cK$ and $e \in \cE$, we conclude that
\be
\label{CcNear}
\sup_{\substack{f \in \cK\\ e \in \cE}}
\|f -R^{\rm near}(L(f)+e)\|_{\cC(D)}
\le (1+2\gamma) \times \mu_{\wt{V},\wt{I}}(\wt{L}) \times \eps ,
\ee
which is the required objective \eqref{ObjNear} with $C_\gamma = 1+2\gamma$.
\epf

\brk
Solving the optimization problem  \eqref{Main} exactly to produce $a^{(j)} \in \bR^m$ may not be possible.
However,
one can solve \eqref{TrunSDP} instead and produce $a^{(N,j)} \in \bR^m$ satisfying
\eqref{Reprod} and a substitute of~\eqref{EqMin} taking the form
\be
\bigg\| Q_j - \sum_{i=1}^m a_i^{(N,j)} \ell_i \bigg\|_{\cC(D)^*} + \f{\eta}{\eps} \|a^{(N,j)}\|_{p'}
\le \mu_{\wt{V},\wt{Q_j}}(\wt{L})
+ \delta^{(N,j)}
\ee
for some quantifiably small quantities $\delta^{(N,j)} \ge 0$,
see Theorem \ref{ThmQuant}.
The linear map \eqref{RNear} with $a^{(N,j)}$ in lieu of $a^{(j)}$ is still be a near-optimal recovery map for the full approximation problem over $\cK$ and~$\cE$.
The previous argument indeed still shows that \eqref{CcNear} holds with $C_\gamma = 1 + 2 \gamma$ loosely
replaced by $C_\gamma  = 1 + 2 \gamma +\max_j \delta^{(N,j)}$.
\erk


\end{document}